\documentclass[11pt]{article}
\usepackage[letterpaper, margin=1.0in]{geometry}
\usepackage{bm}
\usepackage{xcolor}
\usepackage{comment}

\makeatletter

\usepackage[numbers]{natbib}
\usepackage{appendix}

\usepackage{enumerate}
\usepackage[ruled,linesnumbered]{algorithm2e}
\usepackage{multirow}
\usepackage{subfig}
\usepackage{tabularx}
\usepackage{makecell}
\usepackage{booktabs}

\usepackage{amsmath}
\usepackage{amstext}
\usepackage{amsthm}
\usepackage{amssymb}
\usepackage{bbm}
\usepackage{amsbsy}
\usepackage{mathtools}

    \numberwithin{prop}{section} 
    
    \numberwithin{lem}{section} 
    
    \numberwithin{rem}{section} 
    
    \numberwithin{cor}{section} 
    
    \numberwithin{thm}{section} 
    
    \numberwithin{defn}{section}
    
    \numberwithin{fact}{section}

    \theoremstyle{plain}
    \newtheorem*{thm*}{Theorem}
    
    \newcommand{\beq}{\begin{equation}}
    \newcommand{\eeq}{\end{equation}}
    \newcommand{\beqnn}{\begin{equation*}}
    \newcommand{\eeqnn}{\end{equation*}}
        
    \def\bu {\boldsymbol{u}}
    \def\bs {\boldsymbol{s}}
    \def\bv {\boldsymbol{v}}
    
    \def\bx {\boldsymbol{x}}
    \def\by {\boldsymbol{y}}

    \def\R {\mathbb{R}}

    \def\btheta {\boldsymbol{\theta}}
    \def\bthetahat {\hat{\boldsymbol{\theta}}}
    \def\thetahat {\hat{\theta}}

    \newcommand{\norm}[1]{\left\Vert{#1} \right\Vert}

    \newcommand{\matr}[1]{\bm{#1}}

    \def\sign {\mathrm{sign~}}
    
    \DeclareMathAlphabet{\mymathbb}{U}{BOONDOX-ds}{m}{n}
    
    \def\Rn {\R^n}
    \def\Rm {\R^m}
    
    \DeclareMathOperator*{\argmin}{arg\,min}
    \DeclareMathOperator*{\argmax}{arg\,max}

    \newcommand{\normsq}[1]{\left\Vert{#1} \right\Vert_{2}^{2}}
    \newcommand{\normtwo}[1]{\left\Vert{#1} \right\Vert_{2}}
    \newcommand{\normone}[1]{\left\Vert{#1} \right\Vert_{1}}

\usepackage[backref=page]{hyperref}
\renewcommand*{\backref}[1]{}
\renewcommand*{\backrefalt}[4]{%
    \ifcase #1 (Not cited.)%
    \or        (Cited on page~#2.)%
    \else      (Cited on pages~#2.)%
    \fi}

\title{Efficient and robust high-dimensional sparse logistic regression via nonlinear primal-dual hybrid gradient algorithms}
\author{J\'er\^ome Darbon and Gabriel P. Langlois}

\begin{document}

\maketitle

\begin{abstract} 

Logistic regression is a widely used statistical model to describe the relationship between a binary response variable and predictor variables in data sets. It is often used in machine learning to identify important predictor variables. This task, variable selection, typically amounts to fitting a logistic regression model regularized by a convex combination of $\ell_1$ and $\ell_{2}^{2}$ penalties. Since modern big data sets can contain hundreds of thousands to billions of predictor variables, variable selection methods depend on efficient and robust optimization algorithms to perform well. State-of-the-art algorithms for variable selection, however, were not traditionally designed to handle big data sets; they either scale poorly in size or are prone to produce unreliable numerical results. It therefore remains challenging to perform variable selection on big data sets without access to adequate and costly computational resources. In this paper, we propose a nonlinear primal-dual algorithm that addresses these shortcomings. Specifically, we propose an iterative algorithm that provably computes a solution to a logistic regression problem regularized by an elastic net penalty in $O(T(m,n)\log(1/\epsilon))$ operations, where  $\epsilon \in (0,1)$ denotes the tolerance and $T(m,n)$ denotes the number of arithmetic operations required to perform matrix-vector multiplication on a data set with $m$ samples each comprising $n$ features. This result improves on the known complexity bound of $O(\min(m^2n,mn^2)\log(1/\epsilon))$ for first-order optimization methods such as the classic primal-dual hybrid gradient or forward-backward splitting methods. 
\end{abstract}

\begin{centering}
\subsubsection*{Significance statement}
\end{centering}

Logistic regression is a widely used statistical model to describe the relationship between a binary response variable and predictor variables in data sets. With the trends in big data, logistic regression is now commonly applied to data sets whose predictor variables range from hundreds of thousands to billions. State-of-the-art algorithms for fitting logistic regression models, however, were not traditionally designed to handle big data sets; they either scale poorly in size or are prone to produce unreliable numerical results. This paper proposes a nonlinear primal-dual algorithm that provably computes a solution to a logistic regression problem regularized by an elastic net penalty in $O(T(m,n)\log(1/\epsilon))$ operations, where  $\epsilon \in (0,1)$ denotes the tolerance and $T(m,n)$ denotes the number of arithmetic operations required to perform matrix-vector multiplication on a data set with $m$ samples each comprising $n$ features. This result improves on the known complexity bound of $O(\min(m^2n,mn^2)\log(1/\epsilon))$ for first-order optimization methods such as the classic primal-dual hybrid gradient or forward-backward splitting methods.


\bigbreak

\section{Introduction}\label{sec:intro}
Logistic regression is a widely used statistical model to describe the relationship between a binary response variable and predictor variables in data sets~\cite{hosmer2013applied}. It is often used in machine learning to identify important predictor variables~\cite{el2020comparative,zanon2020sparse}. This task, variable selection, typically amounts to fitting a logistic regression model regularized by a convex combination of $\ell_1$ and $\ell_{2}^{2}$ penalties. Variable selection is frequently applied to problems in medicine~\cite{BAGLEY2001979,bursac2008purposeful, greene2014big,PEREIRA2009S199,prive2018efficient,wu2009genome,zhang2018variable}, natural language processing~\cite{berger1996maximum,manning2003optimization,genkin2007large,pranckevivcius2017comparison,taddy2013multinomial}, economics~\cite{lowe2004logistic,theodossiou1998effects,zaghdoudi2013bank,zaidi2016forecasting}, and social science~\cite{achia2010logistic,king2001logistic,muchlinski2016comparing}, among others.

Since modern big data sets can contain up to billions of predictor variables, variable selection methods require efficient and robust optimization algorithms to perform well~\cite{l2017machine}. State-of-the-art algorithms for variable selection methods, however, were not traditionally designed to handle big data sets; they either scale poorly in size~\cite{chu2007map} or are prone to produce unreliable numerical results~\cite{bringmann2018homotopy,loris2008_package,yuan2010comparison,yuan2012improved}. These shortcomings in terms of efficiency and robustness make variable selection methods on big data sets essentially impossible without access to adequate and costly computational resources~\cite{demchenko2013addressing,sculley2014machine}. Further exacerbating this problem is that machine learning applications to big data increasingly rely on computing power to make progress~\cite{dhar2020carbon,kambatla2014trends,l2017machine,leiserson2020there}. Without efficient and robust algorithms to minimize monetary and energy costs, these shortcomings prevent scientific discoveries. Indeed, it is expected that progress will rapidly become economically and environmentally unsustainable as computational requirements become a severe constraint~\cite{thompson2021deep}.

This paper proposes a novel optimization algorithm that addresses the shortcomings of state-of-the-art algorithms used for variable selection. Our proposed algorithm is an accelerated nonlinear variant of the classic primal-dual hybrid gradient (PDHG) algorithm, a first-order optimization method initially developed to solve imaging problems~\cite{esser2010general,pock2009algorithm,zhu2008efficient,chambolle2011first,hohage2014generalization,chambolle2016ergodic}. Our proposed accelerated nonlinear PDHG algorithm, which is based on the work the authors recently provided in~\cite{langlois2021accelerated}, uses the Kullback--Leibler divergence to efficiently fit a logistic regression model regularized by a convex combination of $\ell_1$ and $\ell_{2}^{2}$ penalties. Specifically, our algorithm provably computes a solution to a logistic regression problem regularized by an elastic net penalty in $O(T(m,n)\log(1/\epsilon))$ operations, where  $\epsilon \in (0,1)$ denotes the tolerance and $T(m,n)$ denotes the number of arithmetic operations required to perform matrix-vector multiplication on a data set with $m$ samples each comprising $n$ features. This result improves on the known complexity bound of $O(\min(m^2n,mn^2)\log(1/\epsilon))$ for first-order optimization methods such as the classic primal-dual hybrid gradient or forward-backward splitting methods. 

\subsection*{Organization of this preprint}
In Section~\ref{sec:prelim}, we describe how variable selection works with logistic regression regularized by the elastic net penalty, why this problem is challenging, what the state-of-the-art algorithms are, and what their limitations are. In Section~\ref{sec:methodology}, we describe our approach for solving this problem using the Kullback--Leibler divergence, we derive an explicit algorithm for solving this problem, and we explain why our algorithm overcomes the limitations of current state-of-the-art algorithms. We also describe how our approach can be adapted to solve a broad class logistic regression problems regularized by an appropriate convex penalty, including, for example, the Lasso penalty. Finally, Section~\ref{sec:material} provides a detailed derivation of the explicit algorithms described in Section~\ref{sec:methodology}.


\section{Preliminaries}\label{sec:prelim}
\subsection*{Description of the problem}
Suppose we receive $m$ independent samples $\{(\bx_i,y_i)\}_{i=1}^{m}$, each comprising an $n$-dimensional vector of predictor variables $\bx_i \in \Rn$ and a binary response variable $y_i \in \{0,1\}$. The predictor variables are encoded in an $m \times n$ matrix $\matr A$ whose rows are the vectors $\bx_i = (x_{i1},\dots, x_{in})$, and the binary response variables are encoded in an $m$-dimensional vector $\by$. The goal of variable selection is to identify which of the $n$ predictor variables best describe the $m$ response variables. A common approach to do so is to fit a logistic regression model regularized by a convex combination of $\ell_1$ and $\ell_{2}^{2}$ penalties:
\begin{equation}\label{eq:binary_LR_wp}
    \inf_{\btheta \in \Rn} f(\btheta;\alpha,\lambda) = \inf_{\btheta \in \Rn} \left\{ \frac{1}{m}\sum_{i=1}^{m}\log\left(1 +  \exp{((\matr A \btheta)_i)}\right) - \frac{1}{m}\left\langle\by,\matr A \btheta \right\rangle + \lambda \left(\alpha \normone{\btheta} + \frac{1-\alpha}{2}\normsq{\btheta}\right)\right\},
\end{equation}
where $\lambda > 0$ is a tuning parameter and $\alpha \in (0,1)$ is a fixed hyperparameter. The function $\btheta \mapsto \lambda (\alpha \normone{\btheta} + (1-\alpha)\normsq{\btheta}/2)$ is called the elastic net penalty~\cite{zou2005regularization}. It is a compromise between the ridge penalty ($\alpha = 0$)~\cite{hoerl1970ridge} and the lasso penalty ($\alpha = 1$)~\cite{tibshirani1996regression}. The choice of $\alpha$ depends on the desired prediction model; for variable selection its value is often chosen to be close to but not equal to one~\cite{tay2021elastic}.  

The elastic net regularizes the logistic regression model in three ways. First, it ensures that the logistic regression problem~\eqref{eq:binary_LR_wp} has a unique solution (global minimum)~\cite[Chapter II, Proposition 1.2]{ekeland1999convex}. Second, the $\ell_{2}^{2}$ penalty shrinks the coefficients of correlated predictor variables toward each other (and zero), which alleviates negative correlation effects (e.g., high variance) between highly correlated predictor variables. Third, the $\ell_1$ penalty promotes sparsity in the solution of~\eqref{eq:binary_LR_wp}; that is, the global minimum of~\eqref{eq:binary_LR_wp} has a number of entries that are identically zero~\cite{Foucart2013,el2020comparative,zanon2020sparse}. We note that other penalties are sometimes used in practice to promote sparsity, including, for example, the group lasso penalty~\cite{meier2008group}. In any case, the non-zero entries are identified as the important predictor variables, and the zero entries are discarded. The number of non-zero entries itself depends on the value of the fixed hyperparameter $\alpha$ and the tuning parameter $\lambda$. 

In most applications, the desired value of $\lambda$ proves challenging to estimate. To determine an appropriate value for it, variable selection methods first compute a sequence of minimums $\btheta^{*}(\lambda)$ of problem~\eqref{eq:binary_LR_wp} from a chosen sequence of values of the parameter $\lambda$ and then choose the parameter that gives the preferred minimum~\cite{bringmann2018homotopy,friedman2010regularization}. Variable selection methods differ in how they choose the sequence of parameters $\lambda$ and how they repeatedly compute global minimums of problem~\eqref{eq:binary_LR_wp}, but the procedure is generally the same. The sequence of parameters thus computed is called a regularization path~\cite{friedman2010regularization}.

Unfortunately, computing a regularization path to problem~\eqref{eq:binary_LR_wp} can be prohibitively expensive for big data sets. To see why, fix $\alpha \in (0,1)$ and $\lambda > 0$, and let $\btheta_{\epsilon}(\alpha,\lambda) \in \Rn$ with $\epsilon > 0$ denote an $\epsilon$-approximate solution to the true global minimum $\btheta^{*}(\alpha,\lambda)$ in~\eqref{eq:binary_LR_wp}, i.e.,
\[
f(\btheta_{\epsilon}(\alpha,\lambda);\alpha,\lambda) - f(\btheta^{*}(\alpha,\lambda);\alpha,\lambda) < \epsilon.
\]
Then the best achievable rate of convergence for computing $\btheta_{\epsilon}(\alpha,\lambda)$ in the Nesterov class of optimal first-order methods is linear, that is, $O(\log(1/\epsilon))$ in the number of iterations~\cite{Nesterov2018}. While optimal, this rate of convergence is difficult to achieve in practice because it requires a precise estimate of the largest singular value of the matrix $\matr{A}$, a quantity essentially impossible to compute for large matrices due to its prohibitive computational cost of $O(\min{(m^2 n,mn^2)})$ operations~\cite{hastie2009elements}.  This issue generally makes solving problem~\eqref{eq:binary_LR_wp} difficult and laborious. As computing a regularization path entails repeatedly solving problem~\eqref{eq:binary_LR_wp} for different values of $\lambda$, this process can become particularly time consuming and resource intensive for big data sets.

In summary, variable selection methods work by repeatedly solving an optimization problem that can be prohibitively computationally expensive for big data sets. This issue has driven much research in the development of robust and efficient algorithms to minimize costs and maximize performance.

\subsection*{Algorithms for variable selection methods and their shortcomings} 
The state of the art for computing regularization paths to problem~\eqref{eq:binary_LR_wp} is based on coordinate descent algorithms~\cite{friedman2007pathwise,friedman2010regularization,hastie2021glmnet,simon2011regularization,simon2013blockwise,tibshirani2012strong,wu2008coordinate,yuan2012improved}. These algorithms are implemented, for example, in the popular glmnet software package~\cite{hastie2021glmnet}, which is available in the Python, MATLAB, and R programming languages. Other widely used variable selection methods include those based on the least angle regression algorithm and its variants~\cite{efron2004least,hesterberg2008least,lee2006efficient,tibshirani2013lasso,zou2005regularization}, and those based on the forward-backward splitting algorithm and its variants~\cite{beck2009fast,chambolle2016introduction,daubechies2004iterative,shi2010fast,shi2013linearized}. Here, we focus on these algorithms, but before doing so we wish to stress that many more algorithms have been developed to compute minimums of~\eqref{eq:binary_LR_wp}; see \cite{bertsimas2019sparse,el2020comparative,li2020survey,vidaurre2013survey,zanon2020sparse} for recent surveys and comparisons of different methods and models. 

Coordinate descent algorithms are considered the state of the art because they are scalable, with steps in the algorithms generally having an asymptotic space complexity of at most $O(mn)$ operations. Some coordinate descent algorithms, such as those implemented in the glmnet software~\cite{hastie2021glmnet}, also offer options for parallel computing. Despite these advantages, coordinate descent algorithms generally lack robustness and good convergence properties. For example, the glmnet implementation depends on the sparsity of the matrix $\matr{A}$ to converge fast~\cite{zou2005regularization}, and it is known to be slowed down when the predictor variables are highly correlated~\cite{friedman2007pathwise}. This situation often occurs in practice, and it would be desirable to have a fast algorithm for this case. Another issue is that the glmnet implementation approximates the logarithm term in problem~\eqref{eq:binary_LR_wp} with a quadratic in order to solve the problem efficiently. Without costly step-size optimization, which glmnet avoids to improve performance, the glmnet implementation may not converge~\cite{friedman2010regularization,lee2006efficient}.  Case in point, \citet{yuan2010comparison} provides two numerical experiments in which glmnet does not converge. Although some coordinate descent algorithms recently proposed in~\cite{catalina2018accelerated} and in~\cite{fercoq2016optimization} can provably solve the logistic regression problem~\eqref{eq:binary_LR_wp} (with parameter $\alpha = 1$), in the first case, the convergence rate is strictly less than the achievable rate, and in the second case, the method fails to construct meaningful regularization paths to problem~\eqref{eq:binary_LR_wp}, in addition to having large memory requirements.

The least angle regression algorithm is another popular tool for computing regularization paths to problem~\eqref{eq:binary_LR_wp}. This algorithm, however, scales poorly with the size of data sets because the entire sequence of steps for computing regularization paths has an asymptotic space complexity of at most $O(\min{(m^2 n + m^3,mn^2 + n^3)})$ operations~\cite{efron2004least}. It also lacks robustness because, under certain conditions, it fails to compute meaningful regularization paths to problem~\eqref{eq:binary_LR_wp}~\cite{bringmann2018homotopy,loris2008_package}. Case in point, \citet{bringmann2018homotopy} provides an example for which the least angle regression algorithm fails to converge.

The forward-backward splitting algorithm and its variants are widely used because they are robust and can provably compute $\epsilon$-approximate solutions of~\eqref{eq:binary_LR_wp} in at most $O(\log(1/\epsilon))$ iterations. To achieve this convergence rate, the step size parameter in the algorithm needs to be fine-tuned using a precise estimate of the largest singular value of the matrix $\matr A$. As mentioned before, however, computing this estimate is essentially impossible for large matrices due to its prohibitive computational cost, which has an asymptotic computational complexity of at most $O(\min{(m^2 n,mn^2)})$ operations. Line search methods and other heuristics are often employed to bypass this problem, but they come at the cost of slowing down the convergence of the forward-backward splitting algorithm. Another approach is to compute a crude estimate of the largest singular value of the matrix $\matr{A}$, but doing so dramatically reduces the speed of convergence of the algorithm. This problem makes regularization path construction methods based on the forward-backward splitting algorithm and its variants generally inefficient and impractical for big data sets.

In summary, state-of-the-art and other widely used variable selection methods for computing regularization paths to problem~\eqref{eq:binary_LR_wp} either scale poorly in size or are prone to produce unreliable numerical results. These shortcomings in terms of efficiency and robustness make it challenging to perform variable selection on big data sets without access to adequate and costly computational resources. This paper proposes an efficient and robust optimization algorithm for solving~problem~\eqref{eq:binary_LR_wp} that addresses these shortcomings.
\section{Methodology}\label{sec:methodology}
We consider the problem of solving the logistic regression problem~\eqref{eq:binary_LR_wp} with $\alpha \in (0,1)$. Our approach is to reformulate problem~\eqref{eq:binary_LR_wp} as a saddle-point problem and solve the latter using an appropriate primal-dual algorithm. Based on work the authors recently provided in~\cite{langlois2021accelerated}, we propose to use a nonlinear PDHG algorithm with Bregman divergence terms tailored to the logistic regression model and the elastic net penalty in~\eqref{eq:binary_LR_wp}. Specifically, we propose to use the Bregman divergence generated from the negative sum of $m$ binary entropy functions. This divergence is the function $D_{H}\colon \Rm \times \Rm \to [0,+\infty]$ given by
\begin{equation}\label{eq:kl-divergence}
    D_{H}(\bs,\bs') = \begin{dcases}
     & \sum_{i=1}^{m} s_{i}\log\left(\frac{s_{i}}{s_{i}'}\right) + (1-s_{i})\log\left(\frac{1-s_{i}}{1-s_{i}'}\right) \quad \mathrm{if}\, \bs,\bs' \in [0,1]^m, \\
     & +\infty,\quad \mathrm{otherwise}.
\end{dcases}
\end{equation}
We also show how to adapt our approach for solving the logistic regression problem~\eqref{eq:binary_LR_wp} with the lasso penalty ($\alpha = 0$) and a broad class of convex penalties, including for example the group lasso.

\subsection*{Numerical optimization algorithm}
The starting point of our approach is to express the logistic regression problem~\eqref{eq:binary_LR_wp} in saddle-point form. To do so, we use the convex conjugate formula of the sum of logarithms that appears in~\eqref{eq:binary_LR_wp}, namely
\begin{equation}\label{eq:conv_conj}
\psi(\bs) = \sup_{\bu \in \Rn} \left\{\left\langle \bs,\bu \right\rangle - \sum_{i=1}^{m}\log(1+\exp{(u_i)})\right\} =
\begin{dcases}
& \sum_{i=1}^{m}s_i\log(s_i) + (1-s_i)\log(1-s_i) \quad \mathrm{if}\, \bs \in [0,1]^m, \\
& +\infty,\quad \mathrm{otherwise}.
\end{dcases}
\end{equation}
Hence we have the representation
\[
\sum_{i=1}^{m}\log(1+\exp{((\matr{A}\btheta)_{i})}) = \sup_{\bs \in [0,1]^m} \left\{\left\langle \bs,\matr{A}\btheta \right\rangle - \psi(\bs)\right\},
\]
and from it we can express problem~\eqref{eq:binary_LR_wp} in saddle-point form as
\begin{equation}\label{eq:binary_LR_saddle}
    \inf_{\btheta \in \Rn}\sup_{\bs \in [0,1]^{m}} \left\{ - \frac{1}{m}\psi(\bs) - \frac{1}{m}\left\langle\by - \bs,\matr A \btheta \right\rangle + \lambda \left(\alpha \normone{\btheta} + \frac{1-\alpha}{2}\normsq{\btheta}\right) \right\}.
\end{equation}

A solution to the convex-concave saddle-point problem~\eqref{eq:binary_LR_saddle} is called a saddle point. For $\alpha \in (0,1)$, the saddle-point problem~\eqref{eq:binary_LR_wp} has a unique saddle point $(\btheta^{*},\bs^{*})$, where the element $\btheta^{*}$ itself is the unique global minimum of the original problem~\eqref{eq:binary_LR_wp}~\cite[Proposition 3.1, page 57]{ekeland1999convex}. Hence for our purpose it suffices to compute a solution to the saddle-point problem~\eqref{eq:binary_LR_saddle}, and to do so we can take advantage of the fact that the saddle point $(\btheta^{*},\bs^{*})$ satisfies the following optimality conditions:
\begin{equation}\label{eq:opt_conds}
\frac{1}{m}\matr{A}^T(\by-\bs^{*}) - \lambda(1-\alpha)\btheta^{*} \in \lambda\alpha \partial \normone{\btheta^{*}} \quad \mathrm{and} \quad s_{i}^{*} = \frac{1}{1 + \exp{(-(\matr{A}\theta^{*})_{i})}}\quad \mathrm{for}\,i \in \{1,\dots,m\}.
\end{equation}

The next step of our approach is to split the infimum and supremum problems in~\eqref{eq:binary_LR_saddle} with an appropriate primal-dual scheme. We propose to alternate between a nonlinear proximal ascent step using the Kullback--Leibler divergence~\eqref{eq:kl-divergence} and a proximal descent step using a quadratic function:
\begin{equation}
\begin{alignedat}{1}\label{eq:kl-nPDHG-alg}
     &\bs^{(k+1)} = \argmax_{\bs \in (0,1)^m} \left\{-\psi(\bs) + \left\langle \bs, \matr{A}(\btheta^{(k)} + \rho(\btheta^{(k)} - \btheta^{(k-1)}))\right\rangle - \frac{1}{\sigma}D_{H}(\bs,\bs^{(k)})\right\}, \\
     &\btheta^{(k+1)} = \argmin_{\btheta \in \Rn} \left\{\left(\lambda_{1}\normone{\btheta} + \frac{\lambda_{2}}{2}\normsq{\btheta}\right) + \left\langle \bs^{(k+1)}-\by, \matr{A}\btheta\right\rangle + \frac{1}{2\tau}\normsq{\btheta-\btheta^{(k)}}\right\},
\end{alignedat}
\end{equation}
where $\lambda_{1} = m\lambda\alpha$, $\lambda_{2} = m\lambda(1-\alpha)$, and $\rho,\sigma,\tau > 0$ are parameters to be specified in the next step of our approach. The scheme starts from initial values $\bs^{(0)} \in (0,1)^m$ and $\btheta^{(-1)} = \btheta^{(0)} \in \Rn$. 

The key element in this primal-dual scheme is the choice of the Kullback--Leibler divergence~\eqref{eq:kl-divergence} in the first line of~\eqref{eq:kl-nPDHG-alg}. Its choice is motivated by two facts. First, because it is generated from the sum of $m$ binary entropy that appears \emph{explicitly} in the saddle-point problem~\eqref{eq:binary_LR_saddle} as the function $\psi$ defined in~\eqref{eq:conv_conj}, i.e.,
\begin{equation}\label{eq:kullback_def}
D_{H}(\bs,\bs') = \psi(\bs) - \psi(\bs') -\left\langle \bs - \bs', \nabla \psi(\bs') \right\rangle.
\end{equation}
This fact will make the maximization step in~\eqref{eq:kl-nPDHG-alg} easy to evaluate. Second, because it is strongly convex with respect to the $\ell1$-norm in that
\[
D_{H}(\bs,\bs') \geqslant \frac{1}{2}\normone{\bs-\bs'}^{2}
\]
for every $\bs,\bs' \in [0,1]^{m}$, which is a direct consequence of a fundamental result in information theory known as Pinsker's inequality~\cite{beck2003mirror,csiszar1967information,kemperman1969optimum,kullback1967lower,pinsker1964information}. 

The latter fact, notably, implies that the primal-dual scheme~\eqref{eq:kl-nPDHG-alg} alternates between solving a 1-strongly concave problem over the space $(\Rm,\normone{\cdot})$ and a $\lambda_{2}$-strongly convex problem over the space $(\Rn,\normtwo{\cdot})$. The choice of these spaces is significant, for it induces the matrix norm
\begin{equation}\label{eq:matr_parameter}
\norm{\matr{A}}_{op} = \sup_{\normone{\bs} = 1} \normtwo{\matr{A}^{T}\bs} = \max_{i \in \{1,\dots,m\}} \sqrt{\sum_{j=1}^{n}{A_{ij}^{2}}} = \max_{i \in \{1,\dots,m\}} \normtwo{\bx_{i}},
\end{equation}
which can be computed in \emph{optimal} $\Theta(mn)$ time. This is unlike most first-order optimization methods, such as the forward-backward splitting algorithm, where instead the matrix norm is the largest singular value of the matrix $\matr{A}$, which takes $O(\min{(m^2n,mn^2)})$ operations to compute. This point is \emph{crucial}: the smaller computational cost makes it easy and efficient to estimate all the parameters of the nonlinear PDHG algorithm, which is needed to achieve an optimal rate of convergence.

The last step of our approach is to choose the parameters $\rho$, $\sigma$, and $\tau$ so that the iterations in the primal-dual scheme~\eqref{eq:kl-nPDHG-alg} converge. Based on the analysis of accelerated nonlinear PDHG algorithms the authors recently provided in~\cite[Section 5.4]{langlois2021accelerated}, the choice of parameters
\[
\rho = 1 - \frac{\lambda_{2}}{2\norm{\matr{A}}_{op}^2}\left(\sqrt{1+\frac{4\norm{\matr{A}}_{op}^2}{\lambda_{2}}}-1\right), \quad \sigma = \frac{1-\rho}{\rho}, \quad \mathrm{and} \quad \tau = \frac{(1-\rho)}{\lambda_{2}\rho},
\]
ensure that the iterations converge to the unique saddle point $(\btheta^{*},\bs^{*})$ of problem~\eqref{eq:binary_LR_saddle}. In particular, the rate of convergence is linear in the number of iterations, with
\begin{equation}\label{eq:rate}
    \frac{1}{2}\normsq{\btheta^{*} - \btheta^{(k)}} \leqslant \rho^k\left(\frac{1}{2}\normsq{\btheta^{*} - \btheta^{(0)}} + \frac{1}{\lambda_{2}}D_{H}(\bs^{*},\bs^{(0)})\right).
\end{equation}
This convergence rate is optimal: it is the best achievable rate of convergence in the Nesterov class of optimal first-order methods~\cite{Nesterov2018}.

An important feature of our proposed algorithm is that the minimization steps in~\eqref{eq:kl-nPDHG-alg} can be computed \textit{exactly}. Specifically, with the auxiliary variables $\bu^{(k)} = \matr{A}\btheta^{(k)}$ and $v_{i}^{(k)} = \log\left(s_{i}^{(k)}/(1-s_{i}^{(k)})\right)$ for $i \in \{1,\dots m\}$, the steps in algorithm~\eqref{eq:kl-nPDHG-alg} can be expressed explicitly as follows:
\begin{equation}\label{eq:alg_oPDHG_specific}
\begin{dcases}
    \bv^{(k+1)} &= \frac{1}{1+\sigma}\left(\sigma\bu^{(k)} + \sigma\rho\left(\bu^{(k)} - \bu^{(k-1)}\right) + \bv^{(k)}\right) \\
    s^{(k+1)}_{i} &= \frac{1}{1 + \exp{\left(-v_{i}^{(k+1)}\right)}} \quad \mathrm{for}\,i\in\left\{1,\dots,m\right\}, \\
    \bthetahat^{(k+1)} &= \btheta^{(k)} - \tau \matr{A}^T\left(\bs^{(k+1)} - \by\right)\\
    \theta^{(k+1)}_{j} &= \sign{\thetahat^{(k+1)}_j}\max{\left(0,\frac{\left|\thetahat^{(k+1)}_{j}\right| - \lambda_{1}\tau}{1+\lambda_{2}\tau}\right)} \quad \mathrm{for}\,j\in\{1,\dots,n\}\\
    \bu^{(k+1)} &= \matr{A}\btheta^{(k+1)}.
\end{dcases}
\end{equation}
In addition, from the auxiliary variables and the optimality condition on the right in~\eqref{eq:opt_conds}, we have the limit
\[
\lim_{k \to +\infty} \normtwo{\bu^{(k)} - \bv^{(k)}} = 0,
\]
which can serve as a convergence criterion. We refer to \textit{Material and Methods} for the derivation of algorithm~\eqref{eq:alg_oPDHG_specific} from the iterations in~\eqref{eq:kl-nPDHG-alg}.

Our proposed explicit nonlinear PDHG algorithm~\eqref{eq:alg_oPDHG_specific} offers many advantages in terms of efficiency and robustness. First, the computational bottlenecks in algorithm~\eqref{eq:alg_oPDHG_specific} consist of matrix-vector multiplications and the estimation of the induced matrix $\norm{\matr{A}}_{op}$ given by~\eqref{eq:matr_parameter}. If $T(m,n)$ denotes the number of arithmetic operations required to perform matrix-vector multiplication with the matrix $\matr{A}$, then the asymptotic space complexity for computing the iterations in algorithm~\eqref{eq:alg_oPDHG_specific} as well as the induced matrix $\norm{\matr{A}}_{op}$ is at most $O(T(m,n))$ operations. As mentioned before, this is unlike most first-order optimization methods, such as the forward-backward splitting algorithm, where instead the matrix norm is the largest singular value of the matrix $\matr{A}$, which takes $O(\min{(m^2n,mn^2)})$ operations to compute. This fact is crucial because the smaller computational cost makes it easy and efficient to estimate all the parameters of the nonlinear PDHG algorithm, which is needed to achieve an optimal rate of convergence. 

Another advantage of our algorithm is that it exhibits scalable parallelism because the matrix-vector multiplication operations can be implemented via parallel algorithms. This makes it possible to implement our proposed algorithms in a way that takes advantage of emerging hardware, such as field-programmable gate arrays architectures.

Finally, our algorithm also provably computes an $\epsilon$-approximate solution of~\eqref{eq:binary_LR_wp} in $O(\log(1/\epsilon))$ operations~\cite[Section 5.4]{langlois2021accelerated}. The size of the parameter $\rho$ dictates this linear rate of convergence; it depends on the matrix $\matr{A}$, the tuning parameter $\lambda$, and the hyperparameter $\alpha$. Hence the overall complexity required to compute a global minimum of the elastic net regularized logistic regression problem~\eqref{eq:binary_LR_wp} with tolerance $\epsilon \in (0,1)$ is on the order of $O(T(m,n)\log(1/\epsilon))$ operations.

With these advantages, algorithm~\eqref{eq:alg_oPDHG_specific} overcomes the limitations of the state-of-the-art and other widely-used algorithms for solving the logistic regression problem~\eqref{eq:binary_LR_wp}. We are unaware of any other algorithm that offers these advantages in terms of efficiency and robustness simultaneously.

In general, the nonlinear PDHG algorithm~\eqref{eq:alg_oPDHG_specific} can be adapted to any regularized logistic regression problem for which the penalty is strongly convex on the space $(\Rn,\normtwo{\cdot})$. To do so, substitute this penalty for the elastic net penalty in the minimization problem of the scheme~\eqref{eq:kl-nPDHG-alg} and use its solution in place of the third and fourth lines in the explicit algorithm~\eqref{eq:alg_oPDHG_specific}.

\subsection*{Special case: Logistic regression with the lasso penalty}
In some situations, it may be desirable to fit the regularized logistic regression model~\eqref{eq:binary_LR_wp} without the $\ell{2}^{2}$ penalty ($\alpha = 1$). In this case, algorithm~\eqref{eq:alg_oPDHG_specific} does not apply since it depends on the strong convexity of $\ell{2}^{2}$ penalty. We present here an algorithm for fitting a logistic regression model regularized by an $\ell1$ penalty, or in principle, any convex penalty that is not strongly convex, such as the group lasso.

The $\ell{1}$-regularized logistic regression problem is
\begin{equation}\label{eq:lr_l1}
    \inf_{\btheta \in \Rn} \left\{\frac{1}{m}\sum_{i=1}^{m}\log\left(1 +  \exp{((\matr A \btheta)_i)}\right)  - \frac{1}{m}\left\langle\by,\matr A \btheta \right\rangle + \lambda \normone{\btheta} \right\},
\end{equation}
and its associated saddle-point problem is
\begin{equation}\label{eq:binary_LR_saddle_l1}
    \inf_{\btheta \in \Rn}\sup_{\bs \in [0,1]^{m}} \left\{ - \frac{1}{m}\psi(\bs) - \frac{1}{m}\left\langle\by - \bs,\matr A \btheta \right\rangle + \lambda \normone{\btheta}\right\}.
\end{equation}
The $\ell1$ penalty in~\eqref{eq:lr_l1} guarantees that problem~\eqref{eq:lr_l1} has at least one solution. Accordingly, the saddle-point problem~\eqref{eq:binary_LR_saddle_l1} also has at least one saddle point. As before, we split the infimum and supremum problems in~\eqref{eq:binary_LR_saddle_l1} by alternating between a nonlinear proximal ascent step using the Kullback--Leibler divergence~\eqref{eq:kl-divergence} and a proximal descent step using a quadratic function, but this time we also update the stepsize parameters at each iteration:
\begin{equation}
\begin{alignedat}{1}\label{eq:kl-lasso-alg}
     &\bs^{(k+1)} = \argmax_{\bs \in (0,1)^m} \left\{-\psi(\bs) + \left\langle \bs, \matr{A}(\btheta^{(k)} + \rho^{(k)}(\btheta^{(k)} - \btheta^{(k-1)}))\right\rangle - \frac{1}{\sigma^{(k)}}D_{H}(\bs,\bs^{(k)})\right\}, \\
     &\btheta^{(k+1)} = \argmin_{\btheta \in \Rn} \left\{m\lambda \normone{\btheta} + \left\langle \bs^{(k+1)}-\by, \matr{A}\btheta\right\rangle + \frac{1}{2\tau^{(k)}}\normsq{\btheta-\btheta^{(k)}}\right\}, \\
     &\rho^{(k+1)} = 1/\sqrt{1+\sigma^{(k)}}, \quad \sigma^{(k+1)} = \rho^{(k+1)}\sigma^{(k)}, \quad \tau^{(k+1)} = \tau^{(k)}/\rho^{(k+1)}.
\end{alignedat}
\end{equation}
The scheme starts from initial stepsize parameters $\rho^{(0)} \in (0,1)$, $\tau^{(0)} > 0$ and $\sigma^{(0)} = 1/(\tau^{(0)}\norm{\matr{A}}_{op}^{2})$, and initial values $\bs^{(0)} \in (0,1)^m$ and $\btheta^{(-1)} = \btheta^{(0)} \in \Rn$. 

The following accelerated nonlinear PDHG algorithm computes a global minimum of~\eqref{eq:lr_l1}:
\begin{equation}\label{eq:alg_l1}
\begin{dcases}
    \bv^{(k+1)} &= \frac{1}{1+\sigma^{(k)}}\left(\sigma^{(k)}\bu^{(k)} + \sigma^{(k)}\rho^{(k)}\left(\bu^{(k)} - \bu^{(k-1)}\right) + \bv^{(k)}\right) \\
    s^{(k+1)}_{i} &= \frac{1}{1+\exp{\left(-v_{i}^{(k+1)}\right)}} \quad \mathrm{for}\,i\in\left\{1,\dots,m\right\}, \\
    \bthetahat^{(k+1)} &= \btheta^{(k)} - \tau^{(k)}\matr A^T\left(\bs^{(k+1)} - \by\right)\\
    \theta^{(k+1)}_{j} &= \sign{\thetahat^{(k+1)}_j}\max{\left(0,\left|\thetahat^{(k+1)}_{j}\right| - m\lambda\tau^{(k)}\right)} \quad \mathrm{for}\,j\in\{1,\dots,n\}, \\
    \bu^{(k+1)} &= \matr{A}\btheta^{(k+1)}, \\
    \rho^{(k+1)} &= 1/\sqrt{1+\sigma^{(k)}}, \quad \sigma^{(k+1)} = \rho^{(k+1)}\sigma^{(k)}, \quad \tau^{(k+1)} = \tau^{(k)}/\rho^{(k+1)}.
\end{dcases}
\end{equation}
In addition, from the auxiliary variables and the optimality condition on the right in~\eqref{eq:opt_conds}, we have the limit
\[
\lim_{k \to +\infty} \normtwo{\bu^{(k)} - \bv^{(k)}} = 0,
\]
which can serve as a convergence criterion. The derivation of algorithm~\eqref{eq:alg_l1} from the iterations in~\eqref{eq:kl-lasso-alg} follows from the derivation of algorithm~\eqref{eq:alg_oPDHG_specific} from the iterations in~\eqref{eq:kl-nPDHG-alg} described in $\textit{Material and Methods}$ by setting $\alpha = 1$. According to results provided by the authors in~\cite[Proposition 5.2]{langlois2021accelerated}, the sequence of iterates $\{(\btheta^{(k)},\bs^{(k)})\}_{k=1}^{+\infty}$ converges to a saddle point $(\btheta^{*},\bs^{*})$ of~\eqref{eq:binary_LR_saddle_l1} at a sublinear rate $O(1/k^2)$ in the iterations. Moreover, this sublinear rate satisfies the lower bound
\[
\frac{2\tau^{(0)}\norm{\matr{A}}^{2}_{op}}{1 + 2\tau^{(0)}\norm{\matr{A}}^{2}_{op}}k + \frac{2\tau^{(0)}}{(1 + 2\tau^{(0)}\norm{\matr{A}}^{2}_{op})^2}k^2.
\]
In particular, the constant term multiplying $k^2$ is maximized when $\tau^{(0)} = 1/(2\norm{\matr{A}}_{op}^{2})$. This suggests a practical choice for the free parameter $\tau^{(0)}$.

In general, the nonlinear PDHG algorithm~\eqref{eq:alg_l1} can be adapted to any regularized logistic regression problem for which the penalty is proper, lower semicontinuous and convex, and for which a solution exists. To do so, substitute the $\ell1$ penalty in the minimization problem of the scheme~\eqref{eq:kl-nPDHG-alg} and use its solution in place of the third and fourth lines in the explicit algorithm~\eqref{eq:alg_l1}.
\section{Material and Methods}\label{sec:material}

\subsection*{Derivation of the explicit algorithm~\eqref{eq:alg_oPDHG_specific}}
We derive here the explicit algorithm~\eqref{eq:kl-nPDHG-alg} from the iterations in~\eqref{eq:alg_oPDHG_specific}. Consider the first line of~\eqref{eq:kl-nPDHG-alg}. This maximization problem has a unique maximum inside the interval $(0,1)^m$~\cite[Proposition 3.21-3.23, Theorem 3.24, Corollary 3.25]{bauschke2003bregman}, and the objective function is differentiable. Thus it suffices to compute the gradient with respect to $\bs$ and solve for $\bs$ to compute its global maximum. To do so, it helps to first rearrange the objective function. Substitute $\bs^{(k)}$ for $\bs'$ in equation~\eqref{eq:kl-divergence}, use equation~\eqref{eq:kullback_def}, and rearrange to obtain the objective function
\[
\begin{alignedat}{1}
-\psi(\bs) &+ \left\langle \bs, \matr{A}(\btheta^{(k)} + \rho(\btheta^{(k)} - \btheta^{(k-1)}))\right\rangle - \frac{1}{\sigma}D_{H}(\bs,\bs^{(k)})\\
&= -\psi(\bs) + \left\langle \bs, \matr{A}(\btheta^{(k)} + \rho(\btheta^{(k)} - \btheta^{(k-1)}))\right\rangle - \frac{1}{\sigma}\left(\psi(\bs) - \psi(\bs^{(k)}) - \left\langle \bs-\bs^{(k)},\nabla\psi(\bs^{(k)}) \right\rangle)\right) \\
&= -\left(1 + \frac{1}{\sigma}\right)\psi(\bs) + \left\langle \bs, \matr{A}(\btheta^{(k)} + \rho(\btheta^{(k)} - \btheta^{(k-1)}))\right\rangle + \frac{1}{\sigma}\left\langle \bs-\bs^{(k)},\nabla\psi(\bs^{(k)}) \right\rangle + \psi(\bs^{(k)}).
\end{alignedat}
\]
The optimality condition is then
\[
\nabla\psi(\bs^{(k+1)}) = \frac{\sigma}{1+\sigma}\left(\matr{A}(\btheta^{(k)} + \rho(\btheta^{(k)} - \btheta^{(k-1)}))\right) + \frac{1}{1+\sigma}\nabla\psi(\bs^{(k)}),
\]
where $(\nabla\psi(\bs))_i = \log\left(s_{i}^{(k+1)}/(1-s_{i}^{(k+1)})\right)$ for $i \in \{1,\dots,m\}$ and $\bs \in (0,1)^m$. With the auxiliary variables $\bu^{(k)} = \matr{A}\btheta^{(k)}$ and $v_{i}^{(k)} = \log\left(s_{i}^{(k)}/(1-s_{i}^{(k)})\right)$ for $i \in \{1,\dots m\}$, the optimality condition can be written as
\[
\bv^{(k+1)} = \frac{1}{1+\sigma^{(k)}}\left(\sigma^{(k)}\bu^{(k)} + \sigma^{(k)}\rho^{(k)}\left(\bu^{(k)} - \bu^{(k-1)}\right) + \bv^{(k)}\right).
\]
This gives the first line in~\eqref{eq:alg_oPDHG_specific}. The second line follows upon solving for $\bs^{(k+1)}$ in terms of $\bv^{(k+1)}$. The fifth line follows from the definition of the auxiliary variable $\bu^{(k)}$.

Now, consider the second line of~\eqref{eq:kl-nPDHG-alg}. Complete the square and multiply by $\tau/(1+\lambda_{2}\tau)$ to get the equivalent minimization problem
\[
\btheta^{(k+1)} = \argmin_{\btheta \in \Rn} \left\{\frac{\lambda_{1}\tau}{1 + \lambda_{2}\tau}\normone{\btheta} + \frac{1}{2}\normsq{\btheta - \left(\btheta^{(k)} - \tau\matr{A}^T(\bs^{(k+1)} - \by)\right)/(1+\lambda_{2}\tau)}\right\}.
\]
The unique minimum is computed using the soft thresholding operator~\cite{daubechies2004iterative,figueiredo2001wavelet,lions1979splitting}. With the notation
\[
\bthetahat^{(k+1)} = \btheta^{(k)} - \tau\matr A^T\left(\bs^{(k+1)} - \by\right),
\]
the soft thresholding operator is defined component-wise by
\begin{equation*}
\theta^{(k+1)}_{i} = \sign{\thetahat^{(k+1)}_j}\max{\left(0,\frac{\left|\thetahat^{(k+1)}_{j}\right| - \lambda_{1}\tau}{1+\lambda_{2}\tau}\right)} \quad \mathrm{for}\,j\in\{1,\dots,n\}.
\end{equation*}
The third and fourth lines of~\eqref{eq:alg_oPDHG_specific} are precisely these two equations.




\bibliographystyle{plainnat}
\bibliography{proj-bib}
\end{document}